\documentstyle[12pt,fullpage]{article}

\newtheorem{lem}{Lemme} \newtheorem{th}{Theorem}
\newtheorem{theo}{Th\'{e}or\`{e}me}

\title{Complexes inattendus de droites de saut} \author{Jean Vall\`es}
\begin{document}

\maketitle \begin{abstract} \begin{center} {\bf Unexpected complex of
jumping lines} \end{center} \bigskip We prove here the following results
: \begin{th} Let $E$ a rank 2 vector bundle over ${\bf P}_3$, if $C$ is
a reduced irreducible curve of  ${\bf P}_3^{\vee}$ such that $E_H$ is
unstable for all $H\in C$ then $C$ is a line. \end{th} We define now the
set $W(E)$ as the set of planes $H$ such that the restricted bundle
$E_H$ is unstable (that means non semi-stable). \begin{th} Let $E$  a
rank 2 vector bundle over ${\bf P}_3$, with first chern class
$c_1=c_1(E)$, $L$ a line and   an integer $n\ge 0$. The following
conditions are equivalent : \begin{description} \item[(i)]
$L^{\vee}\subset  W(E)$ and $H^0(E_H(-n+[-c_1/2]))\neq 0$ for a general
point  $H\in L^{\vee}$. \item[(ii)] There exist $m>0$ and a section $t\in
H^0(E(m+[-c_1/2]))$ such that the zero variety of $t$ contains the
 infinitesimal neighbourhood of order $(m+n-1)$ of  $L$.
\end{description}

\end{th} \end{abstract}

\section{Introduction.} Soit $E$ un fibr\'e stable de rang 2 de
premi\`ere classe de chern $c_1\le -3$. Les droites $L$ telles que
$H^0(E_L)\neq 0$ forment un sous sch\'ema de la grassmanienne des
droites de ${\bf P}_3$ dont la codimension attendue est au moins 2 (cf.
\cite{G-P.2}). Pourtant on connait de nombreux exemples de fibr\'es pour
lesquels ce sous sch\'ema est un complexe (voir exemple 1.2).\\ Dans
\cite{G-P} Gruson et Peskine annoncent le r\'esultat suivant :\\ \\{\bf
Lemme B.5 : }{\it Soit $E$ un fibr\'e de rang 2 sur ${\bf P}_3$ tel que
$c_1(E)\le -3$ et  $H^0(E)=0$. Si $K$ est un complexe r\'eduit
irr\'eductible de droites tel que  $H^0(E_L)\neq 0$ pour $L\in K$, il
existe  une courbe r\'eduite irr\'eductible $\Gamma$ de ${\bf P}_3$
telle que les droites de  $K$ sont les droites rencontrant $\Gamma$. De
plus il existe un entier $l$ et une section $s$ de $E(l)$  dont la
vari\'et\'e des z\'eros contient le $(l-1)$-{i\`eme} voisinage
infinit\'esimal de $\Gamma$.}\\ \\ \indent Les auteurs me signalent que
la preuve est incompl\`ete.  Plus  pr\'ecis\'ement la d\'emonstration de
l'irr\'eductibilit\'e g\'eom\'etrique de $K$ au dessus de ${\bf P}_3$
sur laquelle cette preuve repose comporte une erreur.\\ \indent Nous
donnons ici une d\'emonstration de cet \'enonc\'e dans le cas o\`u le
complexe $K$ est form\'e des droites d'une famille de dimension 1 de
plans. Cela peut paraitre surprenant car une telle famille est toujours
g\'eom\'etriquement r\'eductible except\'e dans le cas o\`u elle est
form\'ee des plans contenant une droite. C'est bien entendu le seul cas
possible. En effet, consid\'erons la courbe $C$ de ${\bf P}_3^{\vee}$
naturellement associ\'ee au complexe $K$ (voir lemme 1), on montre :

\begin{theo} Soit $E$ un fibr\'e stable de rang 2 sur ${\bf P}_3$, si
$C$ est une courbe r\'eduite irr\'eductible de  ${\bf P}_3^{\vee}$ telle
que $E_H$ est instable pour tout $H\in C$  alors $C$ est une droite.
\end{theo} {\it Remarque.} On  d\'eduit imm\'ediatement de ce
th\'eor\`eme que  l'ensemble des hyperplans $H$ tels que $E_H$ est
instable (c'est \`a dire non semi-stable) not\'e $W(E)$ est de dimension
au plus 1 dans ${\bf P}_3^{\vee}$ (ce r\'esultat est d\'ej\`a prouv\'e
par Coanda \cite{Coanda}).\\ \\ De plus,  on a l'\'enonc\'e suivant,
directement inspir\'e du lemme B.5 : \begin{theo}
 Soit $E$ un fibr\'e stable de rang 2 sur ${\bf P}_3$ de premi\`ere
classe de chern $c_1=c_1(E)$, $L$ une droite et $n$ un entier $\ge 0$.
Les conditions suivantes sont \'equivalentes: \begin{description}
\item[(i)]
  $L^{\vee}\subset  W(E)$ et $H^0(E_H(-n+[-c_1/2]))\neq 0$ pour $H$ un
point g\'en\'eral de $L^{\vee}$. \item[(ii)] Il existe $m>0$ et une
section $t\in H^0(E(m+[-c_1/2]))$ dont la vari\'et\'e des z\'eros
contient le $(m+n-1)$-{i\`eme} voisinage infinit\'esimal de $L$.
\end{description}  \end{theo} \subsection{Remarques pr\'eliminaires.}
Rappelons qu'un fibr\'e $E$ de rang deux sur  ${\bf P}_3$
  est  stable si et seulement si  $H^0(E([-c_1/2]))=0$,  o\`u $[.]$
d\'esigne la partie enti\`ere et $c_1$ la premi\`ere classe de chern de
$E$.\\ On supposera dans les d\'emonstrations des th\'eor\`emes 1 et 2
que $c_1=0 \,\,\mbox{ou} \, -1$.\\ \\ {\it Remarque 1.} Les
th\'eor\`emes 1 et 2 impliquent le lemme B.5 lorsque le complexe $K$ est
form\'e des droites d'une famille de dimension 1 de plans (voir lemme
1).\\ {\it Remarque 2.} Les th\'eor\`emes 1 et 2 sont vrais sur ${\bf
P}_n,\, n\ge 3$.\\ {\it Remarque 3.} Le lemme B.5 est \'evident lorsque
$E$ n'est pas stable.\\ \\ D\'emontrons la remarque 3 pour $c_1(N)=-3$.
Dans ce cas  $H^0(N(1))\neq 0$ et l'unique section non nulle de $N(1)$
s'annule le long d'une courbe $Z$ : $$ 0 \rightarrow O_{{\bf
P}_3}(-1)\rightarrow N\rightarrow {\cal I}_Z(-2)\rightarrow 0.$$ On en
d\'eduit que $L\in K$ si et seulement si $L\cap Z \neq \emptyset$. Le
complexe $K$ \'etant r\'eduit irr\'eductible, ceci prouve l'existence
d'une courbe r\'eduite irr\'eductible $\Gamma$ telle que ${\cal I}_Z
\subset {\cal I}_{\Gamma}.$\\ \\ Le lemme suivant  \'etablit la
correspondance entre le complexe $K$ de droites de saut d'un fibr\'e
$E$, complexe form\'e des droites d'une famille de dimension 1 de plans,
et la courbe $C$ de  $W(E)$. \begin{lem} Soit $N$ un fibr\'e de rang 2
stable tel que $c_1(N)\le -2$, on a l'\'equivalence : $$ H^0(N_L)\neq 0,
\,\,{\rm\mbox{pour tout}}\,\,L\subset H \Longleftrightarrow H^0(N_H)\neq
0.$$ \end{lem} {\it d\'emonstration :} L'implication ($\Leftarrow$)
r\'esulte de la suite exacte suivante : $$0\longrightarrow
N_H(-1)\longrightarrow N_H\longrightarrow N_L\longrightarrow 0.$$
Inversement, d'apr\`es Grauert-Mulich le fibr\'e $N_H$ est instable. Si
$H^0(N_H)=0$, on peut donc supposer que $H^0(N_H(1))\neq 0$, soit,
$$0\longrightarrow O_H(-1)\longrightarrow N_H\longrightarrow {\cal
I}_{Z_H}(c_1+1)\longrightarrow 0,$$ o\`u $Z_H$ est un groupe de points
de $H$. On en d\'eduit que pour une droite g\'en\'erale $L$ de $H$,
$H^0(N_L)=0$, ce qui contredit l'hypoth\`ese.\\ \\ {\small {\bf
Remerciements :} Je remercie Christian Peskine qui a  dirig\'e ce
travail et Iustin Coanda pour ses conseils et les nombreuses et utiles
discussions que nous avons eues.} \subsection{Exemple.} Avant de
d\'emontrer les th\'eor\`emes 1 et 2 donnons un exemple de fibr\'e
${\cal E}$ poss\`edant une famille de dimension 1 de plans instables et
d\'ecrivons   une section de ce fibr\'e ayant les propri\'et\'es
annonc\'ees dans le th\'eor\`eme 2.\\ \\  {\it Exemple.} On reprend
l'exemple de Gruson-Peskine d\'ecrit par Mei-Chu Chang \cite{M.C.C}:
Soient $(r+1)$ plans $H_1,...,H_{r+1}\, (r\geq 2)$
 articul\'es autour d'une droite $D$, et dans chaque plan $H_i$ une
courbe $X_i$ de degr\'e $(2r-1)$;  les $(r+1)$ courbes sont choisies 2
\`a 2 disjointes. On note:
 $$ X=\bigcup_{i=1}^{(r+1)} X_i \quad \mbox{o\`u}\quad X_i\subset H_i .$$
 Il existe un fibr\'e ${\cal E}$ stable de rang 2 tel que $c_1({\cal
E})=0$ et une extension  \begin{equation}
 0\longrightarrow O_{{\bf P}_3}\longrightarrow {\cal E}(r)
\longrightarrow {\cal I}_X(2r)\longrightarrow 0.  \end{equation} En
effet comme $w_X=O_X(2r- 4)$, on a $\mbox{Ext}^1({\cal I}_X,w_{{}_{{\bf
P}_3}}(4-2r))\neq 0$,  et ${\cal E}$ est stable car $H^0({\cal
I}_X(r))=0$.\\
 Soit $H$ un plan g\'en\'eral contenant $D$, alors  $H$ coupe proprement
$X$, et  $(1)$  reste exacte apr\`es tensorisation par $O_H(1-2r)$ \[
 0\longrightarrow O_H(1-2r)\longrightarrow {\cal
E}_H(1-r)\longrightarrow {\cal I}_{X\cap H}(1)\longrightarrow 0. \]
Comme $X\cap H\subset D$ on a  $H^0({\cal I}_{X\cap H}(1))\neq 0$, donc
  $H^0({\cal E}_H(1-r))\neq 0$ (on remarque de plus que $H^0({\cal
E}_H(-r))=0$ pour tout plan $H$).
 Nous avons montr\'e que  $D^{\vee}$ est une droite de ${\bf
P}_3^{\vee}$ dont les points correspondent \`a des plans  instables du
fibr\'e ${\cal E}$. \\ \\ \indent D\'ecrivons maintenant une section de
ce fibr\'e ayant les propri\'et\'es annonc\'ees. La suite exacte (1)
montre que
  $\, H^0({\cal I}_X(r+1))\simeq H^0({\cal E}(1)).$
 La courbe $X$ est  contenue dans la r\'eunion de $(r+1)$ plans. Cette
surface induit une section non nulle $t$ de ${\cal E}(1)$. On note
$\Gamma$ le lieu des z\'eros de la section $t$ (c'est une courbe car
$H^0({\cal E})=0$).\\ Montrons que $\Gamma$ contient le $(r-1)$-{i\`eme}
voisinage infinit\'esimal de $D$. Sinon soit
 $H$ un plan contenant  $D$ et coupant proprement $\Gamma$. Dans ce cas
la section $t_H\in H^0({\cal E}_H(1))$ s'annule en codimension 2 le long
du groupe de points $\Gamma\cap H$ \begin{equation} 0 \longrightarrow
O_H\longrightarrow {\cal E}_H(1)\longrightarrow {\cal I}_{\Gamma\cap H}
(2)\longrightarrow 0. \end{equation} Il en r\'esulte  que  $H^0({\cal
E}_H(-1))=0$ , ce qui contredit $H\in W({\cal E})$. Donc la section
$t_H$ s'annule le long d'une courbe (d'\'equation $f_H=0$) contenant $D$
 et \'eventuellement d'un groupe de points. Si le degr\'e de cette
courbe est $k\le (r-1)$, la section  $t_H/f_H\in H^0({\cal E}_H(1-k))$
s'annule en codimension 2, et on trouve $H^0({\cal E}_H(-k))=0$.
 C'est impossible car $H^0({\cal E}_H(1-r))\neq 0$.\\
  Par cons\'equent la restriction de $\Gamma$ \`a un plan g\'en\'eral
contenant $D$ est une courbe plane de degr\'e $r$. On en d\'eduit que
$\Gamma$ contient le $(r-1)$-{i\`eme} voisinage infinit\'esimal de $D$.\\
\\ Montrons de plus que  le support de $\Gamma$ est contenu dans $D$.
 En effet supposons que $\Gamma$ poss\`ede une autre composante
irr\'eductible $\Gamma_1$.
 Soit $H$ le plan contenant $\Gamma_1$ et $L$. Le lieu des z\'eros de la
section $t_H \in H^0({\cal E}_H(1))$ contient une courbe de degr\'e $\ge
(r+1)$.
 Alors $H^0({\cal E}_H(-r))\neq 0$, ce qui est impossible.\\ \\ Le lieu
des z\'eros de la  section $t\in H^0({\cal E}(1))$ est une courbe
$\Gamma$ telle que : \begin{description} \item{-} $\deg (\Gamma)=r^2+r$
et  $p_a(\Gamma)=1-r^2-r$ (o\`u $p_a$ est le genre arithm\'etique).
\item{-} $\Gamma$ contient le $(r-1)$-i\`eme voisinage infinit\'esimal
de $D$. \item{-} Le support de $\Gamma$ est inclus dans $D$.
\end{description} \section{D\'emonstration du th\'eor\`eme 1}
 La d\'emonstration est divis\'ee en deux \'etapes.
 On v\'erifie tout d'abord que $C$ est plane.  On d\'emontre ensuite que
c'est une droite en interpr\'etant la singularit\'e apparaissant au
point correspondant au plan de ${\bf P}_3^{\vee}$ contenant la courbe.\\
D'apr\`es le th\'eor\`eme de semi continuit\'e la fonction
 $h^0(E_H(-k))$ est minimale sur un ouvert non vide $U$ de $C$. Soit
$n\ge 0$ tel que $h^0(E_H(-n))=1$ pour  un point   $H\in U$
(remarquons que $n\ge 1$ lorsque $c_1(E)=0$). \paragraph*{Etape 1 : } La
d\'emonstration de cette premi\`ere partie repose essentiellement sur la
remarque bien connue suivante que nous ne red\'emontrons pas.\\ \\ {\it
Remarque 1.}  Soit $H$ un plan g\'en\'eral de $C$  et $s_H$ l'unique
section non nulle de $H^0(E_H(-n))$. On note $Z(s_H)$ le sch\'ema des
z\'eros de $s_H$, si $L$  est une droite de $H$ on a : $$\deg (O_{L\cap
Z(s_H)})\ge  r \Longleftrightarrow H^0(E_L(-n-r))\neq 0.$$ \indent On
suppose que $C$ n'est pas plane, i.e. la famille de plans correspondant
aux points de $C$ ne poss\`ede pas de point fixe. Soit $H_i$ et $H_j$
deux points g\'en\'eraux de $U$. Comme $Z(s_{H_i})\cap
Z(s_{H_j})=\emptyset$, on a d'apr\`es la remarque 1, $h^0(E_{H_i\cap
H_j}(-n))=1$. Par cons\'equent les sections $s_{H_i}\in
H^0(E_{H_i}(-n))$ et $s_{H_j}\in H^0(E_{H_j}(-n))$ uniques \`a une
constante pr\`es coincident sur $H_i\cap H_j$ i.e. $s_{H_i}(x)=\lambda
s_{H_j}(x)$ pour $x\in H_i\cap H_j$ et $\lambda \in {\bf C}$.\\ \\ Soit
${\bf F} \subset {\bf P}_3 \times {\bf P}_3^{\vee}$ la vari\'et\'e
d'incidence
 points-plans de ${\bf P}_3$. Posons  $X:= {\bf F}\cap({\bf P}_3 \times
C) $ et consid\'erons
 les projections  : $  X\stackrel{p}\rightarrow{\bf P}_3  $ et $
X\stackrel{q}\rightarrow C $.  Les fibres de $X$ au dessus de $C $ sont
des ${\bf P}_2$, et au dessus de ${\bf P}_3$ elles correspondent aux
sections hyperplanes de $C $. On a : $$h^0(p^*E(-n)_{\mid
q^{-1}(H)})=h^0(E_H(-n))=1 \quad\mbox{pour tout}\,\,H\in U.$$  Le
faisceau $q_*p^*E(-n)$ est  coh\'erent de rang 1 sur $C$ sans torsion
donc inversible. Le morphisme canonique $q^*q_*p^*E(-n)\rightarrow
p^*E(-n)$ s'annule le long d'un sous sch\'ema ferm\'e $Z$ de $X$ de
codimension $\ge 2$. Il induit un morphisme compos\'e $\phi $ au dessus
de ${\bf P}_3$ : $$ X\setminus Z \stackrel{\phi}\longrightarrow {\bf
P}_{{}_{{\bf P}_3}} (E) $$ $$ p\searrow \hspace{0.5cm}\swarrow \pi
\hspace{0.2cm}$$ $${\bf P}_3 $$ L'image de $\phi $ est un ouvert dans un
diviseur r\'eduit irr\'eductible $D$ de  ${\bf P}_{{}_{{\bf P}_3}}(E) $.
Le diviseur $D$ correspond \`a une section  $t\in H^0(S_kE(m))$ (o\`u
$S_kE(m)$ est la $k$-i\`eme puissance sym\'etrique de $E$). La fibre
g\'en\'erale $p^{-1}(x)$ de $ X\setminus Z$ au dessus de ${\bf P}_3$
consiste en $d=\deg C$ points $(x,H_1),...,(x,H_d)$. L'image de
$(x,H_i)$ par $\phi$ est donn\'ee par le point $s_{H_i}(x)$ de ${\bf
P}_{{}_{{\bf P}_3}}(E_x) $ o\`u $s_{H_i}$ est l'unique section non nulle
de $E_{H_i}$. Il r\'esulte de la remarque 1 que les points
$(x,H_1),...,(x,H_d)$ s'envoient sur le m\^eme point  de  ${\bf
P}_{{}_{{\bf P}_3}}(E_x) $. \\ \\On en d\'eduit que la fibre
g\'en\'erale de $D$ au dessus de  ${\bf P}_3$ est irr\'eductible et
qu'elle consiste en un point simple par lissit\'e g\'en\'erique. Le
diviseur $D$ est donc birationnel \`a ${\bf P}_3$. Il en r\'esulte que
$t$ est une section de $E(m)$ avec $m>0$ car $E$ est stable. Quitte \`a
enlever la composante de codimension 1,
 on peut supposer que $t$ s'annule en codimension 2 le long d'une courbe
$\Gamma$.\\ \\Les sections $s_H\in H^0(E_H(-n))$ et $t_H \in
H^0(E_H(m))$ sont proportionnelles. Par hypoth\`ese $
H^0(E_H(-n-1))=0$,  la section $s_H$ s'annule alors en codimension $\ge
2$. Il existe $f_H\in H^0(O_H(m+1))$ tel que $t_H=f_Hs_H$. Mais pour un
plan g\'en\'eral H,   $t_H$ s'annule aussi en codimension 2 (car la
courbe $C$ ne poss\`ede pas de point fixe), ce qui est impossible car
$m>0$.  \paragraph*{Etape 2:} Soit $X$ le plan de ${\bf P}_3^\vee$
contenant la courbe $C$,  et $x\in {\bf P}_3$ le point correspondant \`a
$X$. On consid\`ere l'\'eclatement ${\bf \widetilde{P}}_3 $ de ${\bf
P}_3$ en $x$ et les morphismes naturels $p:{\bf \widetilde{P}}_3
\rightarrow {\bf P}_3 $ et $q:{\bf \widetilde{P}}_3 \rightarrow {\bf
P}_2 $. On a $ {\bf \widetilde{P}}_3=\mbox{Proj}_{{}_{{\bf
P}_3}}(\bigoplus_i  {\cal M}_x^i) $ et $
 {\bf \widetilde{P}}_3={\bf P}_{{}_{{\bf P}_2}}(O_{{\bf P}_2}(1) \oplus
O_{{\bf P}_2} ) $.\\
  On appelle $\Delta$ le diviseur exceptionnel, le morphisme $q_{\mid
\Delta}: \Delta \rightarrow  {\bf P}_2 $ est un isomorphisme. De plus si
$l$ est un point de ${\bf P}_2$ on a : $$ (p^*E)_{\mid q^{-1}(l)} \simeq
E_L. $$
 Le morphisme  $q$ \'etant plat,  le th\'eor\`eme de semi-continuit\'e
implique  que  $ h^0((p^*E)_{\mid q^{-1}(l)})$ atteint son minimum sur
un ouvert non vide $U$ de $ {\bf P}_2$. \noindent Alors  il existe $
a\geq (n+1)$ tel que $ h^0((p^*E)_{\mid q^{-1}(l)}(-a))=h^0(E_L(-a))=1$
pour $l\in U$. L'entier $a$ est $\ge (n+1)$ d'apr\`es la  remarque 1.
\noindent Le faisceau $q_*p^*E(-a) $ est un faisceau coh\'erent de rang
1, il est r\'eflexif donc inversible sur $ {\bf P}_2$. Soit  $k$ tel
que  $q_*p^*E(-a)=O_{{\bf P}_2}(-k)$.
 On consid\`ere le morphisme canonique:
 $$\widetilde{s}:q^*O_{{\bf P}_2}(-k)=q^*q_*p^*E(-a)\longrightarrow
p^*E(-a).$$ Il est injectif. En effet pour $l\in U$, $ p^*E(-a)_{\mid
q^{-1}(l)} \simeq O_L\bigoplus O_L(-a+c_1)$. D'apr\`es le th\'eor\`eme
de changement de base, la restriction de $\widetilde{s}$ \`a la fibre de
$l$ est en fait le morphisme d'\'evaluation (de l'unique section non
nulle de $E_L(-a)$) qui ne s'annule pas,
          $$\widetilde{s}_{\mid q^{-1}(l)}:H^0(E_L(-a))\otimes
O_L\longrightarrow E_L(-a).$$  Le lieu des z\'eros du  morphisme
$\widetilde{s}$ est alors contenu dans la r\'eunion des fibres
$q^{-1}(l)$ o\`u $L$ est une droite dont l'ordre de saut est $>a$. On
sait (cf lemme 9, \cite{Barth}) que ce lieu des z\'eros ne contient pas
d'hypersurface, c'est \`a dire qu'il ne contient au plus qu'un nombre
fini de fibres. \\ \\ Le morphisme $\widetilde{s}$  induit  donc une
section non nulle $\widetilde{t}$ de  $ p^*E(-a)\otimes q^*O_{{\bf
P}_2}(k)$ qui s'annule le long d'une courbe de $ {\bf \widetilde{P}}_3$.
   On remarque que  $H^0(E(-a))= 0$ implique $k>0$.\\
 Comme $p_*q^*O_{{\bf P}_2}(k) \simeq {\cal M}_x^k(k)$, l'image directe
sur ${\bf P}_3$ de la section $\widetilde{t}$ est
 une section non nulle $t\in H^0(E(-a)\otimes {\cal M}_x^k(k))$, ce qui
d\'emontre $k>a$. La section $t$ induit   une suite exacte: $$
0\longrightarrow O_{{\bf P}_3}\stackrel{t} \longrightarrow E(k-a)
                        \longrightarrow {\cal
I}_{\Gamma}(2k-2a+c_1)\longrightarrow 0 $$ o\`u $\Gamma$ est une courbe
qui contient le $(k-1)$-i\`eme voisinage
 infinit\'esimal de $x$, i.e ${\cal I}_{\Gamma}\subset {\cal M}_x^k$.\\
Soient $H$ un point g\'en\'eral de $C$, $t_H$ la restriction de $t$ au
plan $H$ et $s_H$ l'unique section non nulle de $E_H(-n)$. On
consid\`ere le morphisme de fibr\'es sur $H$ suivant: $$ O_H(n)\oplus
O_H(a-k)\stackrel{(s_H,t_H)} \longrightarrow E_H. $$
 Le d\'eterminant de ce morphisme  est une section de ${\cal
M}_{x,H}^k(c_1+k-a-n)$. Etant donn\'e que $(c_1 +k-a-n)<k$ cette section
est nulle i.e $s_H \wedge t_H \equiv 0$ \\ Le lieu des z\'eros de la
section $s_H$ \'etant de codimension deux, on en d\'eduit qu'il existe
$$ f_H\in H^0(O_H(k-a+n))\ /\ t_H=f_H s_H.$$ Ceci montre que la courbe
du plan $H$ d'equation $ f_H=0 $ est contenue dans $\Gamma$. En
cons\'equence tout plan de $C$ contient une composante irr\'eductible
de  $\Gamma$. Il est clair que  $\Gamma$ contient une droite $L$ telle
que  $C=L^{\vee}$.  \section{D\'emonstration du th\'eor\`eme 2} {\bf
(ii)} $\Rightarrow$ {\bf (i)}
 On note $t_H\in H^0(E_H(m))$ la restriction de $t$  \`a un plan
g\'en\'eral $H$ contenant $L$. Par hypoth\`ese $\,{\cal I}_{\Gamma}.O_H
\subset {\cal I}_{L,H}^{(m+n)}\simeq O_H(-m-n) $. En d'autres termes,
la vari\'et\'e des z\'eros de $t_H$ contient une hypersurface
(d'\'equation $f_H=0$) de degr\'e $(m+n)$. Alors $t_H/f_H$ est une
section non nulle de $H^0(E_H(-n))$.\\ \\ {\bf (i)} $\Rightarrow$ {\bf
(ii)}
 Dans le cas pr\'esent, le complexe de droites associ\'e \`a  $L^{\vee}$
(droites rencontrant $L$)   est g\'eom\'etriquement irr\'eductible. Nous
reprenons donc les id\'ees de la d\'emonstration du lemme B.5 pour
\'etablir le r\'esultat.  \\ \\ On reprend la construction de l'\'etape
1 (th\'eor\`eme 1). Consid\'erons la vari\'et\'e d'incidence
 points-plans ${\bf F}$ de ${\bf P}_3$. Posons  $X:= {\bf F}\cap({\bf
P}_3 \times  L^{\vee}) $ et consid\'erons
 les projections  : $  X\stackrel{p}\rightarrow{\bf P}_3  $ et $
X\stackrel{q}\rightarrow L^{\vee} $. Comme pr\'ec\'edemment on a un
morphisme compos\'e $\phi $ au dessus de ${\bf P}_3$ : $$ X\setminus Z
\stackrel{\phi}\longrightarrow {\bf P}_{{}_{{\bf P}_3}} (E),$$ et
l'image de $ X\setminus Z$ est un diviseur $D$ de ${\bf P}_{{}_{{\bf
P}_3}} (E) $.
 Le complexe \'etant g\'eom\'etriquement irr\'eductible  $D$  est
birationnel \`a ${\bf P}_3$.
 Il en r\'esulte que $D$ correspond \`a une section $t$  de $E(m)$ avec
$m>0$ car $E$ est stable. La section $t$ s'annule le long d'une courbe
$\Gamma$.\\ \\ Soit $H\supset L$, les sections $s_H\in H^0(E_H(-n))$ et
$t_H \in H^0(E_H(m))$ sont proportionnelles. Par hypoth\`ese $
H^0(E_H(-n-1))=0$,  la section $s_H$ s'annule alors en codimension $\ge
2$. On en d\'eduit qu'il existe $f_H\in H^0(O_H(m+n))$ telle que
$t_H=f_Hs_H$.\\ \\ Ceci montre que la courbe du plan $H$ d'\'equation
$f_H=0$ est contenue dans $\Gamma$. Par cons\'equent  tout plan $H$ de
$L^{\vee}$ contient une composante irr\'eductible de $\Gamma$ de degr\'e
$(m+n)$. Il est clair que la courbe  $\Gamma$ contient le
$(m+n-1)$-i\`eme voisinage infinit\'esimal de $L$.

\end{document}